\documentclass[a4paper, 10pt, oneside]{article}
\usepackage[margin=1in]{geometry}
\usepackage{times}
\usepackage{amsfonts,amsmath,amssymb,amsthm}
\usepackage{graphicx} 
\newtheorem{df}{Definition}
\newtheorem{teo}{Theorem}
\newtheorem{ex}{Example}
\newtheorem{lem}{Lemma}

\title{A comparison of arithmetical operations with $f$ correlated fuzzy numbers}
\author{Diogo Sampaio da Silva \& Roberto Antonio Cordeiro Prata }
\date{\today}

\begin{document}

\maketitle

\begin{abstract}
We present a brief introduction to a class of interactive fuzzy numbers, called $f$-correlated fuzzy numbers, which consist of pairs of fuzzy numbers where one is dependent on the other by a continuous monotone injective function. We have deduced some equations that can directly calculate the results of the sums and products of $f$-correlated fuzzy numbers, using only basic operations with real numbers, intervals on the real line and the function that relates the fuzzy numbers being considered. We proved that their correlated and standard sum coincide, and that in a certain sense, the correlated product is contained in the standard product.

\end{abstract}

\section{Introduction}
Fuzzy Mathematics is a recent area in mathematics when compared to already established mathematical fields because it emerged in the second half of the 20th century, from the work of \cite{Zadeh1965}. The initial idea was to expand the notion of set, traditionally defined by the dichotomy between membership and non-membership. Thus, to formally deal with propositions and uncertain quantities, the concept of fuzzy set was introduced. A fuzzy set, as opposed to a classical or crisp set, is a set with continuous degrees of membership, modeling real situations where there is no exact criterion delineating a class of objects \cite{Barros2017,Zadeh1965}. Therefore, a fuzzy set is represented by a function in the form $\varphi_{A} : U \to [0,1]$. Here, the domain is a universal set and the function's output represents the degree of membership, increasing progressively.

More specifically, the notion of a fuzzy number was developed to generalize the concept of a real number in a similar manner. If $U$ is a topological space, we can define the $\alpha$-levels of a fuzzy subset $A$, for $\alpha \in [0,1]$. If $\alpha >0$, then $[A]^{\alpha}$ is the set of elements $x \in U$ with $\varphi_{A} (x) \ge \alpha$, and if $\alpha =0$, then $[A]^{0}$ is the closure of the set of elements $x \in U$ with $\varphi_{A} (x) > 0$. A fuzzy number is a fuzzy subset of $\mathbb{R}$ such that all its $\alpha$-levels are non-empty compact intervals. An essential property is that a fuzzy number can be identified as a collection of embeddable compact intervals, indexed in the interval $[0,1]$. By applying Zadeh's Extension Principle, we derive operations among fuzzy numbers, utilizing closed intervals of real numbers \cite{Barros2017}. 

However, such operations can also be derived through possibility distributions \cite{1375791}. Even when discussing only one among these operations, the literature already proposes various methods to define, for instance, the differences between fuzzy numbers, such as: traditional, via distributions, Hukuhara's, generalized Hukuhara's, generalized, and CIA \cite{de2015diferenccas}. An important use of these operations, specially between interactive fuzzy numbers is the study of interactive derivative \cite{DEBARROS201764}.

Through the joint possibility distribution of two fuzzy numbers, we explore the so-called interactive fuzzy numbers, which can be described as those whose values are interdependently assigned \cite{fuller2004interactive}.

Literature already exists on how to add interactive fuzzy numbers, especially those that are completely correlated \cite{coroianu2013additivity}, as well as studies on the multiplication of interactive fuzzy numbers \cite{coroianu2013multiplication}. Our aim is to further complement such studies in a more general field, through a specific yet comprehensive case of interactive fuzzy numbers: $f$-correlated fuzzy numbers.

An important tool is the extension principle for interactive fuzzy numbers, presented as follows:
\begin{equation}
f_{C}(A_{1}, \dots , A_{n}) (y) = \sup_{y=f(x_{1}, \dots , x_{n})} C(x_{1}, \dots , x_{n}).
\label{extensão}
\end{equation}
Here, $C$ is a joint possibility distribution with marginal possibility distributions being fuzzy numbers $A_{1}, \dots , A_{n}$, and $f$ is a continuous function from the $n$-dimensional Euclidean space to the real line.

The concept of $f$-correlated fuzzy numbers is a generalization of the concept of linearly correlated fuzzy numbers, now using a monotonic injective function instead of a linear function. 

\begin{df}[$f$-correlated fuzzy numbers \cite{fcorrelacionado, CabralPrataBarros}]
Let \( f : X \longrightarrow Y \), \( X, Y \subset \mathbb{R} \) be an injective monotone function and \( A, B  \) fuzzy numbers. We say that \( A \) and \( B \) are \emph{correlated via function} or \emph{\( f \)-correlated}, if their joint possibility distribution \( C \) is given by
\begin{equation}
\varphi_C(x, y) = \varphi_A(x) \chi_{\{y=f(x)\}}(x, y) = \varphi_B(y) \chi_{\{y=f(x)\}}(x, y)
\end{equation}
where,
\begin{equation}
    \chi_{\{y=f(x)\}}(x, y) = 
\begin{cases} 
1 & \text{if } y = f(x) \\
0 & \text{if } y \neq f(x) 
\end{cases}.
\end{equation}
\end{df}

Using the extension principle for interactive fuzzy numbers, we have our main basic results.

\begin{lem}[\cite{fcorrelacionado, CabralPrataBarros}]
\label{lem1}
If $A$ and $B$ are $f$-correlated fuzzy numbers, then for each $\alpha \in [0,1]$, we have $[B]^\alpha = f([A]^\alpha)$.
\end{lem}

\begin{lem}[\cite{fcorrelacionado, CabralPrataBarros}]
\label{lem2}
Let $A,B$ be $f$-correlated fuzzy numbers with joint possibility distribution $C$, and let $g(x,y)=x+y$, $h(x,y)=xy$. Then the $\alpha$-cuts of their correlated sum and product are given by:
\begin{align}
[A +_f B]^\alpha &= \overline{\{x + f(x) ; \varphi_A(x) > \alpha\}} \\
[A \cdot_f B]^\alpha &= \overline{\{x f(x) ; \varphi_A(x) > \alpha\}}
\end{align}
where $\overline{X}$ is the closure of $X$, and the membership functions satisfy:
\begin{equation}
\varphi_{A +_c B}(z) = \sup_{z=x+y} \varphi_C(x,y) \quad \textrm{and} \quad \varphi_{A \cdot_c B}(z) = \sup_{z=xy} \varphi_C(x,y)
\end{equation}
with $\varphi_C(x,y) = \varphi_A(x)\chi_{\{y=f(x)\}}(x,y)$.
\end{lem}

\section{Results and Discussions}
Hereinafter, $A$ and $B$ are $f$ correlated fuzzy numbers, and we study comparisons between their standard and correlated sums and products. Also, we represent the $\alpha$-levels of a $A$ by $[A]^\alpha =[L(\alpha), U(\alpha)]$. 

Finally, we observe that continuous monotone injective functions have the following property: $f([a,b]) = [f(a), f(b)]$, if $f$ is increasing, and $f([a,b]) = [f(b), f(a)]$ , if $f$ is decreasing. That is not hard to prove: we know that the image of an interval by a continuous function is an interval; and by Weierstrass extreme value theorem, we also know that the image of a bounded and closed interval by a continuous function is other bounded and closed interval; finally, if $a \le x \le b$, we conclude by definition of injective function.

\subsection{Standard sums and products of $f$ correlated fuzzy numbers}

As a consequence of the Lemma \ref{lem1}, for all $\alpha \in [0,1]$, we have:
\begin{align}
[A + B]^\alpha &= [A]^\alpha  + f([A]^\alpha)\\
[A \cdot B]^\alpha &= [A]^\alpha \cdot f([A]^\alpha)
\end{align}

\subsubsection{Correlation via a linear function}
We obtain that, if $f(x)= qx+r$, then, for all $\alpha \in [0,1]$,
\begin{equation}
[A + B]^\alpha = [L(\alpha), U(\alpha)] + [qL(\alpha)+r, qU(\alpha)+r] = [(q+1)L(\alpha)+r, (q+1)U(\alpha)+r] = (q+1) [A]^\alpha +r
\end{equation}
if $f$ is increasing, i.e., $q>0$. However, if $q<0$, i.e., $f$ is decreasing, then:
\begin{equation}
[A + B]^\alpha = [L(\alpha), U(\alpha)] + [qU(\alpha)+r, qL(\alpha)+r] = [(L(\alpha) + qU(\alpha)+r, U(\alpha) +qL(\alpha) + r] = [A]^\alpha - q [A]^\alpha +r
\end{equation}
in any case, we also have:
\begin{equation}
[A \cdot B]^\alpha = [\min P, \max P] 
\end{equation}
where $P=\{qL(\alpha)^2 + rL(\alpha), q U(\alpha) L(\alpha) + r L(\alpha), q U(\alpha) L(\alpha) + rU(\alpha), qU(\alpha)^2 + rU(\alpha) \}$.

\subsubsection{Correlation via a hyperbolic function}
Similarly, we obtain that for $f(x) = q/x + r$, given $\alpha \in [0,1]$, if $f$ is increasing, then
\begin{equation}
[A + B]^\alpha = [L(\alpha), U(\alpha)] + [q/L(\alpha)+r, q/U(\alpha)+r] = [q/L(\alpha) + L(\alpha) +r, U(\alpha) + q/U(\alpha)+r]
\end{equation}
if $f$ is decreasing, then
\begin{equation}
[A + B]^\alpha = [L(\alpha), U(\alpha)] + [q/U(\alpha)+r, q/L(\alpha)+r] =  [q/U(\alpha) + L(\alpha) + r, q/L(\alpha)+ U(\alpha) + r]
\end{equation}
in any case, we obtain
\begin{equation}
[A \cdot B]^\alpha = [\min P, \max P] 
\end{equation}
where $P=\{ q + r L(\alpha), U(\alpha)q/L(\alpha) + r U(\alpha), L(\alpha)q/U(\alpha) + r L(\alpha), q + r U(\alpha)\}$.
\subsection{Some $f$ correlated sums and products}

As a consequence of the Lemma \ref{lem2}, we know that, if $f(x)= qx+r$, then for all $\alpha \in [0,1]$, we obtain:
\begin{align}
[A +_f B]^\alpha &= \overline{\{(q + 1)x+r ; \varphi_A(x) > \alpha\}} = (q+1)[A]^\alpha + r\\
[A \cdot_f B]^\alpha &= \overline{\{x (qx+r) ; \varphi_A(x) > \alpha\}} = q \overline{\{ x^2; \varphi_A(x) > \alpha\}}+r \overline{\{ x; \varphi_A(x) > \alpha\}} = q [A \cdot_{id} A]^\alpha + r[A]^\alpha
\end{align}
where $id$ is the identity function.

Also, if $f(x)=q/x + r$, then for all $\alpha \in [0,1]$, we have:

\begin{align}
[A +_f B]^\alpha &= \overline{\{x + q/x + r ; \varphi_A(x) > \alpha\}} = [A]^\alpha + q \overline{\{1/x ; \varphi_A(x) > \alpha\}} + r  \\
[A \cdot_f B]^\alpha &= \overline{\{x (q/x+r) ; \varphi_A(x) > \alpha\}} =  \overline{\{ q+rx; \varphi_A(x) > \alpha\}} = q + r[A]^\alpha
\end{align}

\subsection{Main Results}

We remember that the infimum and supremum preserve addition, they also preserve multiplication of positive numbers, moreover, the supremum is monotone increasing and the infimum is monotone decreasing, as can be seen in classic real analysis texts \cite{Berberian1994}.

\begin{teo}
The correlated and standard sum of $f$-correlated fuzzy numbers coincide. Also, the $\alpha$-levels of the correlated product are subsets of the standard product. 
\end{teo}

\begin{proof}
Given $A$ and $B$ $f$-correlated fuzzy numbers, for each $\alpha \in [0,1]$, by Lemma \ref{lem2}, and the definition of fuzzy numbers, we have
\begin{align}
[A +_f B]^\alpha & = \overline{\{x + f(x) ; \varphi_A(x) > \alpha\}} \\
&= \left[ \inf_{\varphi_A(x) > \alpha} x + f(x), \sup_{\varphi_A(x) > \alpha} x + f(x)\right]\\
& = \left[ \inf_{\varphi_A(x) > \alpha} x + \inf_{\varphi_A(x) > \alpha} f(x), \sup_{\varphi_A(x) > \alpha} x + \sup_{\varphi_A(x) > \alpha} xf(x)\right]
\\ &= [A]^\alpha + f ([A]^\alpha)
\end{align}

We also have
\begin{equation}
[A \cdot_f B]^\alpha = \overline{\{xf(x) ; \varphi_A(x) > \alpha\}} = \left[ \inf_{\varphi_A(x) > \alpha} x  f(x), \sup_{\varphi_A(x) > \alpha} x f(x)\right]    
\end{equation}
however, since 
\begin{equation}
\{x  f(x)
; \varphi_A(x) > \alpha\} \subset \{x_1  f(x_2); \varphi_A(x_1), \varphi_A(x_2) > \alpha\}    
\end{equation}
we can conclude
\begin{equation}
\left[ \inf_{\varphi_A(x) > \alpha} x  f(x), \sup_{\varphi_A(x) > \alpha} x f(x)\right] \subset \left[ \inf_{\varphi_A(x_1),\varphi_A(x_2) > \alpha} x_1  f(x_2), \sup_{\varphi_A(x_1),\varphi_A(x_2) > \alpha} x_1  f(x_2) \right]      
\end{equation}
that means $[A \cdot_f B]^\alpha \subset [A]^\alpha \cdot f([A]^\alpha)$.
\end{proof}

Now we present a counterexample of the equality of the correlated and standard products. 

\begin{ex}
Given a fuzzy number $A$, with $\text{supp } A = (-2,1)$, we have $[A \cdot A]^0 = [-2, 4]$ and $[A \cdot_{id} A]^0 = [0,2]$.   
\end{ex}

Finally, we observe an interesting property of the correlated sum: the additive inverse. Similarly, there is a multiplicative inverse for the correlated product.

\begin{teo}
Given a fuzzy number $A$, there are fuzzy numbers $B$ and $C$, and monotone injective functions $f$ and $g$, such that $A$ and $B$ are $f$-correlated and $[A +_f B] =0$, and $A$ and $C$ are $g$-correlated and $[A \cdot_g C] =1$   
\end{teo}

\begin{proof}
If $f(x)=-x$, then given $\alpha \in [0,1]$, by Lemma \ref{lem2}, we have
\begin{equation}
[A +_f B]^\alpha  = \overline{\{x -x ; \varphi_A(x) > \alpha\}} = [0,0]   
\end{equation}
similarly, if $g(x)=1/x$, then we obtain
\begin{equation}
[A \cdot_g C]^\alpha  = \overline{\{x/x  ; \varphi_A(x) > \alpha\}} = [1,1]   
\end{equation}
that means the correlated sums and products have inverse elements.
\end{proof}
\section{Conclusions}

Since the distributivity of the sum of real numbers over the usual product fuzzy product does not always hold, from our theorem, we can infer that at least it holds in the case of the sum of fuzzy numbers correlated via a linear function. Our example shows the difference between $[A \cdot_{id} A]$ and $[A \cdot A]$: the $\alpha$-levels of the first are subsets of the $\alpha$-levels of the last.

\section*{Acknowledgments}

This study was financed in part by the Coordenação de Aperfeiçoamento de Pessoal de Nível Superior – Brasil (CAPES) – Finance Code 001.

\end{document}